\documentclass[a4paper]{article}
\usepackage{subfigure}
\usepackage[english]{babel}
\usepackage{amsmath}
\usepackage{graphicx}
\usepackage[colorinlistoftodos]{todonotes}
\usepackage{psfrag,epsf}

\usepackage{amssymb}
\usepackage{color}
\usepackage{subfig}
\usepackage{bm}
\usepackage{booktabs}
\usepackage{graphicx} 
	
\usepackage{color}
\usepackage{soul} 
\usepackage{amsmath}
\usepackage{enumitem}
\usepackage{bigstrut}

\newcommand{\edge}{%
  \mathrel{-}
  \joinrel\joinrel 
  \mathrel{-}
}

\newcommand{\M}{\mbox{$\mathcal{M}$}}

\newcommand{\E}{\mbox{$\mathbb{E}$}}
\newcommand{\I}{\mathbbm{1}}
\newcommand{\R}{\mbox{$\mathbb{R}$}}

\usepackage{url} 
\usepackage{amsthm,mathrsfs,amsfonts}	
\usepackage{enumitem}
\usepackage{bbm}
\newcommand*{\QEDB}{\hfill\ensuremath{\square}}%

\def\defeq{\mathrel{\mathop:}=}

\numberwithin{equation}{section}

\newtheorem{definition}{Definition}[section]
\newtheorem{property}{Property}

\newtheorem{remark}{Remark}

\def\defeq{\mathrel{\mathop:}=}

\usepackage{graphicx} 
\usepackage{subfig}
\usepackage{float}		
\usepackage{color}
\usepackage[normalem]{ulem}
\title{\bf Sensitivity indices for output on a Riemannian  manifold}

\author{Ricardo Fraiman \thanks{rfraiman@cmat.edu.uy}\hspace{.2cm}\\
    Centro de Matem\'atica, Facultad de Ciencias, \\
     Universidad de la Rep\'ublica, Uruguay.\\
   Fabrice Gamboa \\
    Institut de Math\'ematiques de Toulouse, France.\\
 and \\ 
Leonardo Moreno \\ 
Departamento de M\'etodos Cuantitativos, FCEA, \\
Universidad de la Rep\'ublica, Uruguay.}

\date{2018-10-26}

\begin{document}
\maketitle

\begin{abstract}
In the context of computer code experiments, sensitivity analysis of a complicated input-output system is often performed by ranking the so-called Sobol indices. One reason of the popularity of Sobol's approach relies on the simplicity of the statistical estimation of these indices using the so-called Pick and Freeze method. 
In this work we propose and study  sensitivity indices for the case where the  output  lies on a Riemannian manifold. These indices are based on a Cram\'er von Mises like criterion that takes into account the geometry of the output support. We propose a Pick-Freeze like estimator of these indices based on an $U$--statistic. The asymptotic properties of these estimators are studied. Further, we provide and discuss some interesting numerical examples. 

\end{abstract}

\noindent%

\section{Introduction}	

In many situations occurring in applied mathematics (for example mathematical models or  numerical simulation), when working with an input-output system with uncertain (random)  inputs, it is crucial to understand the global influence of one of the inputs on the output. This problem is generally called global sensitivity analysis (or in short sensitivity analysis). We refer, for example to
\cite{tarantola2008} and  \cite{saltelli2000} for an overview on practical aspects of sensitivity analysis.  Sensitivity analysis aims to give some quantitative indicator allowing to score the global influence of each input variable of the model. A very popular index well tailored for the case of a real valued output  is the so-called Sobol index early proposed in \cite{sobol1993}.   It is based on second order moments of the distributions related to the input-output system.
Different strategies have been implemented for its statistical estimation (see for instance \cite{gamboa2013_s} or \cite{rugama2018}).  This kind of second order indices may also be considered in the more general  frame where the output is a  real vector or a real function (see for instance \cite{klein2014} or \cite{marrel2011}).  
Nevertheless, the variance-based indices have some drawbacks. Indeed, they only study the impact on the variance of the output and so provide  only restricted summary of the output distribution (see \cite{daveiga2015}). Hence, 
considering higher order indices involving the whole distributions and not only their second order properties may give more accurate information on the system. 
In \cite{pianosi2015},  a sensitivity measure based on the Kolmogorov-Smirnov statistic is proposed and studied while entropy-based sensitivity measures are considered for instance  
in \cite{liu2006r}  and \cite{borgonovo2007}. 
For a general overview we also refer to  \cite{borgonovo2017}.
When the output $Z=F(X)$ ($X=(X_i)$ is the input vector), is valued in a more general space, the authors of  \cite{borgonovo2016} have proposed the following general sensitivity index for the input $X_i$:
\begin{displaymath}
S_i \defeq \E_{X_i} \left[ d \left(P_Z, P_{Z \vert X_i} \right) \right].
\end{displaymath} 
Here,  $d(\cdot, \cdot)$ is a given dissimilarity measure between  probability measures, $\E_{X_i}(\cdot)$ denotes the expectation with respect to $X_i$ and $P_Z$  (resp. $P_{Z \vert X_i}$) is the unconditional (resp. conditional) probability. Notice that in general, the study of an accurate statistical estimation of such index is not obvious.

Up to our knowledge, the case where the output takes its values on a Riemannian manifold has not yet been studied. 
Beside, in \cite{gamboa2015} a  sensitivity measure in the case of a scalar output and based on the Cram\'er von Mises distance is proposed an studied. This paper extends this last approach to the more general frame of an output valued on a Riemannian manifold.

In the last decades starting from the pioneer 1945' work of Rao \cite{rao1945}, the statistical theory for data valued on a Riemannian manifold have received a lot of interest and contributions. One of the reason of the development of such theory is the spectacular increase of the computation power allowing the treatment of more and more complex objects and structures on computers. References on the subject are numerous. We refer to \cite{mardia1972}, \cite{bhattacharya2012}  and  \cite{ellingson2015} and the references therein for an overview. In this paper, our aim is to bridge this theory to sensitivity analysis. We build a general sensitivity index. On the one hand this index takes into account the whole characteristic of the involved distributions and not only their second order moments.  On the other hand the working frame is a general system with an output valued on a Riemannian manifold. One of the main ingredient of our work is a generalization of the Cram\'er von Mises criterion replacing, before integrating, the half lines by geodesic balls. Consequently, our new index has the nice property to be invariant with respect to the isometric maps of the Riemannian manifold. Moreover, we show that this generalized index can be easily estimated using a Pick Freeze like estimator. This estimator involves a $U$-statistics. Without loss of generality,  we will assume that the random input variables are real and independent. Sensitivity analysis for dependent inputs may be performed but is generally less readable (see for example  \cite{gamboa2015b} and references therein).

 The paper is organized as follows. In Section  \ref{conceptosbasicos}, in order to be self contained, we recall some key facts and properties on Riemannian manifolds. In Section \ref{indicedef}, we define our new sensitivity index for  Riemannian manifold. Further, we will build a Pick an Freeze estimator based on a $U$-statistic and study its asymptotic (consistency), and non asymptotic (concentration inequality) properties.    In Section \ref{datosreales}, we study the sensitivity of a real life model coming from mechanic. The output of the system is the rigidness matrix. This matrix is involved to model the linear elasticity of a solid object.  The set of these matrices is as a sub-manifold of the symmetric positive matrices. In the last section, we bold the main advantages of our method and discuss its possible extensions. The codes in  Julia and R languages (see \cite{julia2012} and \cite{team2017} respectively) are available upon request to the authors. All proofs are postponed to the Appendix.

\section{Basic concepts}

\label{conceptosbasicos}
Let us begin with  some useful tools and facts in Riemannian geometry.  A \textit{Riemannian metric} $g$ on a manifold $\M$ allows to define for every point $p \in \M$  a scalar product $g_p(\cdot,\cdot)$ acting on the tangent space (at $p$), $T_p \M$. This scalar product depends smoothly  on $p$.   The \textit{Riemannian manifold} $(\M,g)$ is the manifold equipped with the Riemannian metric  $g$. 
For any $v \in T_p \M$ the Riemannian norm of $v$ is given by $\Vert v \Vert \defeq \sqrt{g_p(v,v)}$.
Let $x,y\in M$ and  $\gamma: I:=[0, 1] \rightarrow \M$ be a continuously differentiable curve contained in the Riemannian  manifold $\M$ (that is assumed to be connected) with 
$\gamma(0)=x$ and $\gamma(1)=y$. Let denote by  $L(\gamma)$ its length. We may now define the \textit{induced distance}  $d_g(x, y)$ between the points $x$ and $y$ of $\M$  by setting
    $$d_g(x,y) = \inf \{ L(\gamma) : \gamma \, 
 \}.$$
$d_g$ is the Riemannian distance on $\M$  with respect to the $g$ metric. A \textit{geodesic} (with speed $s \in \mathbb{R}^{+}_0$) is a smooth map $\alpha:I \rightarrow \M$, such that $\Vert \alpha^{'}(t) \Vert= s$  for all $t \in I$ and which is {\it locally length minimizing}. For $p \in \M$ and  $v\in T_p \M$, there exists a unique geodesic  $\alpha_{(p,v)}(t)$ starting from that point with initial tangent vector $v$. The \textit{exponential map} is the map $\exp_p$ given by $\exp_p (v) := \alpha_{(p,v)}(1)$.

\begin{figure}[!ht]
\centering
\subfigure{\includegraphics[width=50mm]{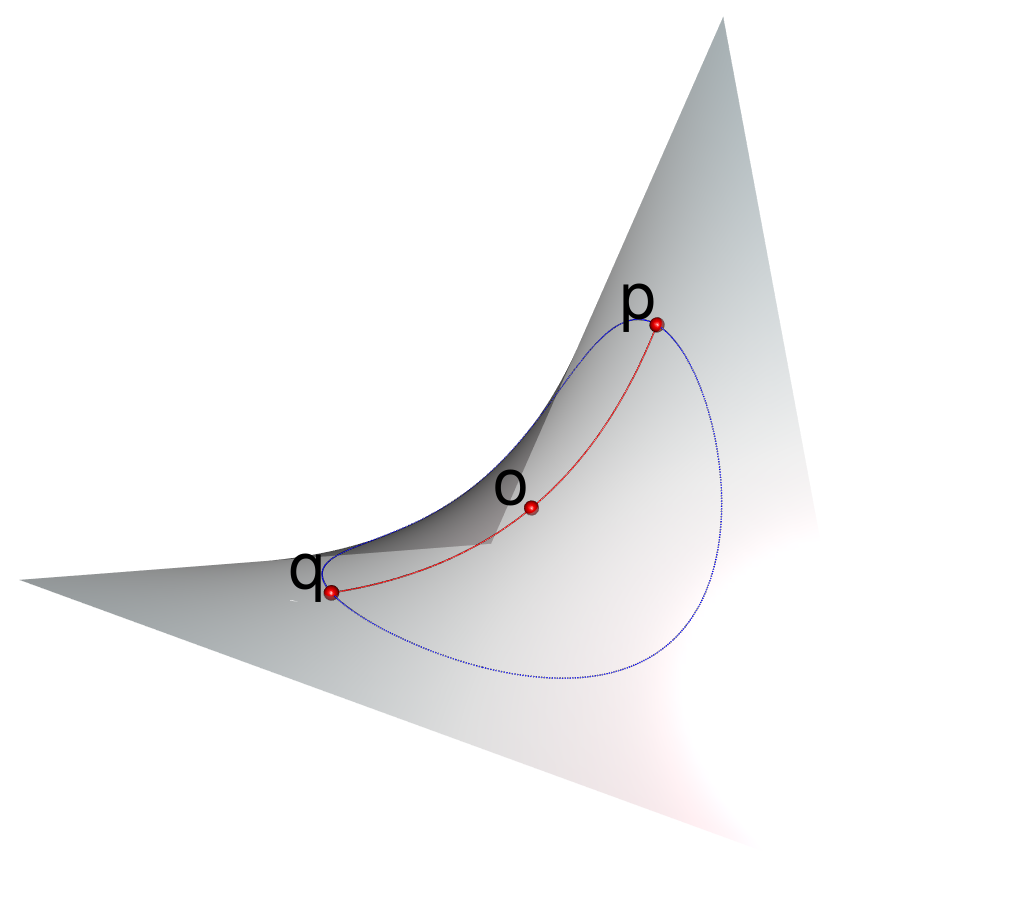}}
\subfigure{\includegraphics[width=40mm]{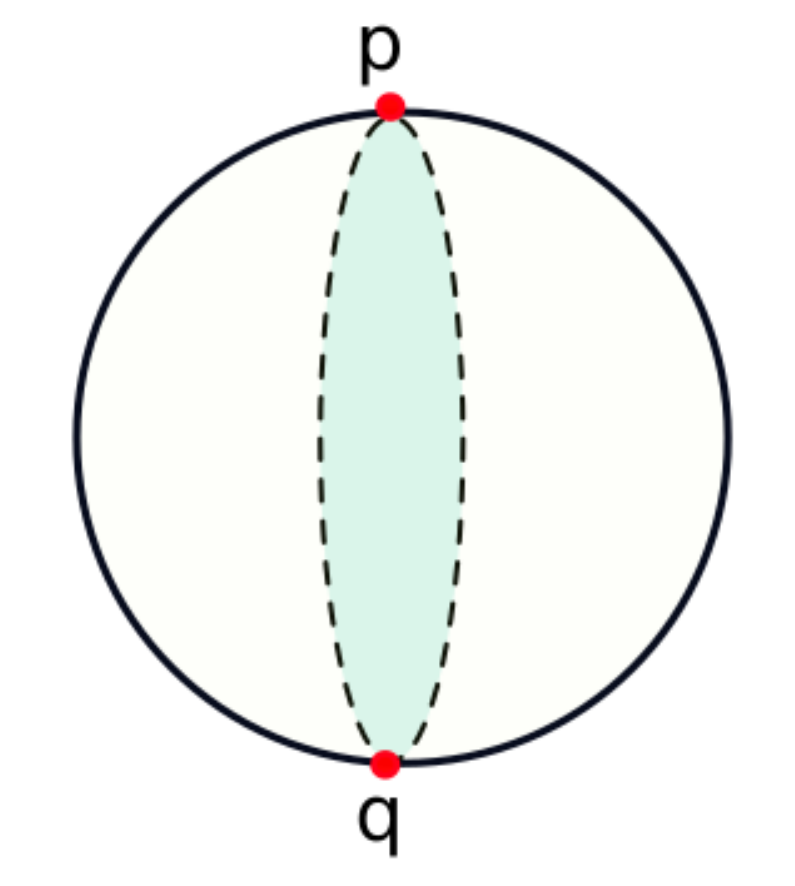}}
\caption{\textit{Left Panel}: Representation of a geodesic ball of diameter $\overline{PQ}$ and center $O$ on a hyperbolic paraboloid. \textit{Right Panel}: A graph of two minimizing geodesics that join two antipodal points of a sphere.}
\label{bola_diametro}
\end{figure}

Notice further that the geodesic that joins two points is not necessarily unique.
The \textit{cut locus}  of $p$ in the tangent space is defined to be the set of all vectors $v \in T_{p} \M$ such that  $\exp _{p}(tv)$ is a minimizing geodesic for  $0 \leq t \leq 1$  but fails to be one for $0 \leq t \leq 1+ \epsilon$, $\epsilon >0$. The \textit{cut locus of $p$ in $\M$}, denoted  by $C_{\M}(p)$, is defined as the image of the cut locus of  $p$ in the tangent space under the exponential map at $p$. The \textit{injectivity radius of $p$} is the maximal radius of centered balls on which the exponential map is a diffeomorphism. The \textit{injectivity radius of the manifold}  $r_{iny}$ is the infimum of the injectivity over the manifold. For example, in the sphere $S_{d}$ of $\R^{d+1}$ the cut locus of a point $p$ is its antipodal point  $-p$. There are infinite minimizing geodesics that connect a point with its antipodal (Figure  \ref{bola_diametro} shows two geodesics) but this configuration has zero probability of being drawn when dealing with the uniform probability measure. 

\begin{figure}[!ht]
\centering
\subfigure{\includegraphics[width=38mm]{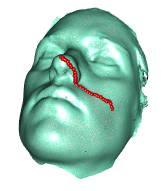}}
\subfigure{\includegraphics[width=50mm]{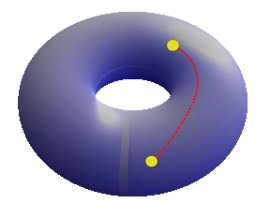}}
\caption{\textit{Left Panel}: A graph of a geodesic that determines two points on the representation
of a face. \textit{Right Panel}: A graph of a geodesic that determines two points on the representation
of a Torus.}
\label{geodesica}
\end{figure}

Let $(\M,g)$ be a connected and orientable Riemannian manifold,(see \cite{docarmo1992}, page $18$). We will assume that  $(\M,d_g)$ is a complete  separable metric space. Since $(\M,d_g)$ is complete, the Hopf--Rinow theorem (see \cite{docarmo1992}, p.146) implies that for any pair of points $p,q \in \M$  there exist at least one geodesic path in $\M$ connecting $p$ and $q$. 
If the manifold $\M$ fulfills the assumptions of the Cartan--Hadamard Theorem (see \cite{petersen2006}, p.162)--that is, if $\M$ is a simply connected, complete Riemannian manifold with non--positive curvature (a Hadamard manifold), then the geodesic is unique.  
 Let $X$ be a random element taking values in $\M$, with distribution $P$. In many useful examples (eg. the sphere),  it occurs that geodesic uniqueness fails.  But in what follows, we will only need that this failure not occurs too often. Roughly speaking, we will assume a condition on $P$ ensuring that there is a unique geodesic between each pair of points $p,q \in \M$ with  probability one. More precisely, we will assume that the random element $Y$ has a density $f_Y$ with respect to the volume measure $d \nu(y)$ on $\M$, (which exists since $\M$ is orientable - see  the last section of \cite{folland2013}) fulfilling the following condition:
Given $q \in \M$ let 
\begin{equation}
A_q\defeq \left\{y \in \M  \Big/
\begin{array}{ll}
\textrm{ there are more than one different } \\ 
\textrm{minimizing geodesic connecting $y$ and $q$}
\end{array} \right \}.
\end{equation}
Then, there exist a Borel set $B_q \subset \M$ with $A_q \subset B_q$ and such that for any  $q \in \M$
\begin{equation}
\label{unicidadd}
\int_{B_q} f_Y(y) d \nu(y)=0.
\end{equation}
\begin{remark}
Let $C_{\M}(p)$ stand for the  cut locus of $p$ in a complete manifold $\M$, (see \cite{docarmo1992}, p. $267$).  If $P \left( Y \in C_{\M}(p)\right) = 0$,  for all   $p \in \M$, then condition  (\ref{unicidadd}) is fulfilled (see \cite{pennec2006}). For instance, if the Riemannian manifold is Hadamard, then the cut locus of any point $p$ is empty and the condition is obviously fulfilled. In the unit sphere $S_{d} \defeq \{u \in \mathbb{R}^d / \Vert u \Vert= 1 \}$  the cut locus of a point $p$ is its opposite point  $-p$ and any  probability measure having a density supported by the sphere fulfills condition  (\ref{unicidadd}).
\end{remark}
In this paper,  we will assume that the Riemannian manifold $(\M,g)$ with the induced distance $d_g$ is connected and oriented, and that the metric space $(\M,d_g)$ is separable and complete. We will also assume  that given two points $p,q \in \M$  there is a unique geodesic determined by $p,q$ with probability one with respect to the tensorial probability measure $d\xi(p,q) \defeq f_Y(p)f_Y(q)d \nu(p)d\nu(q)$, (see Figure \ref{geodesica}).

\subsection{Geodesic balls with diameter $\overline{pq}$}
\label{bolasgeodesicas}
For any pair $p,q \in \M$ that determine a unique geodesic   $\overline{pq}$,  we define the ball of diameter  $\overline{pq}$  as the closed ball whose center is the middle point of the geodesic joining $p$ and $q$ with radius $d_g(p,q)/2$. It will be denoted by $B_{pq}$. Figure  \ref{bola_diametro}  depicts  a geodesic ball of diameter $p$ and $q$, where $o$ is the midpoint of the geodesic.
In \cite{christensen1970} it is shown that the family  $\mathcal B_p:= \{{B}_{p}, p \in \M\}$ of geodesic balls with center $p$ is a determining class; that is, if two probability measures $\eta, \nu$ on $\M$  coincide on any event of this class, then $\eta=\nu$.
We assume that $\M$ is a compact Riemannian manifold. We will start by proving that the family of balls $$\mathcal{B}_{pq}:= \{B_{pq}, p,q \in \M\ \textrm{and there is a unique geodesic between $p$ and $q$} \}$$ is also a determining class if the radius of injectivity  $r_{iny}$ of $\M$ is positive.  We will need the following weaker assumption:

\textbf{HB) B-continuity} A probability measure $\nu$ defined on a Riemannian manifold $\M$ fulfils \bf HB \rm if $\nu(\partial A )=0$  for all closed Borel sets  $A$ on $\M$ (Here $\partial A$ denotes the topological boundary of the set $A$).

\begin{property}
\label{propiedad1}
Let $(\M,g)$ be a compact Riemannian manifold and $\nu$ be a probability distribution fulfilling \textbf{HB}, such that $\nu \left( \partial \M \right)=0$. If the injectivity radius is positive  then  $\mathcal{B}_{pq}$ is a determining class of $\nu$.
\end{property}
 In  \cite{fraiman2017} we showed what if $K \in \M$ is a compact set, then the family of balls 
$\mathcal{B}_{pq}$ is also a Glivenko--Cantelli class in $K$; that is, $$\sup_{p,q \in K}  \vert  P \left(B_{pq} \right) - P_n\left(B_{pq} \right) \vert 
 \rightarrow 0  \quad \textrm{a.s. \ \  as  $n \rightarrow + \infty$}.$$

In the next section, we define a sensitivity index using this family of sets.

\section{Sensitivity Index in Geodesic Balls}
\label{indicedef}
In order to build a sensitivity index, we use an idea previously developed in \cite{gamboa2015}. In our frame we  compare the distribution $P_Y$ and the conditional distribution $P_{Y \vert X_i}$ on the family of geodesic balls on the manifold (instead of half--lines in \cite{gamboa2015}. So that, we build an intrinsic index related to the manifold. This index depends only on the dimension and structure of the manifold and not on the dimension of the space where it is immersed, unlike in \cite{klein2014} and \cite{gamboa2015}.

\subsection{Constructing the index}

Let $\mathbb{X}=(X_1, \ldots, X_d) \in \R^d$ be an random vector and assume that $\mathbb{P} \defeq  \mathbb{P}_1 \times \ldots \times \mathbb{P}_d$ is the probability law of $\mathbb{X}$. Further, let
\begin{displaymath}
Z=f(X_1, \ldots, X_d),
\end{displaymath}
\noindent where $f: \mathbb{R}^d \rightarrow  \M$ is a continuous function. Here, $\M$ a Riemannian manifold of dimension $k$. We wish to understand how sensitive is the  output $Z$ to perturbations in some of the input variables $(X_1,\ldots, X_d)$. As any manifold can be immersed in $\mathbb{R}^p$ for a sufficiently large $p$ (see \cite{nash1954}), this problem could be considered as a particular case of a multivariate output. Nevertheless,  in such a way the geometry and \textit{``minor dimensionality''} of the output are not taken into account. We begin by studying a very particular case and then extend our results to any Riemannian manifold, that satisfies   the previous conditions. Let $h$  a measurable function such that $h(X_1,X_2)$ is integrable, we set $$\E_{X_2}\left(h(X_1,X_2) \right) \defeq \E \left( h(X_1,X_2)  / X_1 \right).$$

\subsection{A very particular case: The real line.}

Let us first focus on the very particular case  $\M=\R$. Let $F$ be the distribution function of $Z$, $F(t) \defeq P(Z \leq t) \, (t\in \mathbb{R})$, and $F^{\nu}$ be  the distribution of $Z$ conditioned by a subset of the variable $(X_1, \ldots, X_d)$. That is, let $\nu= \{i_1,\ldots, i_k\}  \subset \{1, \dots, d\}$, and

\begin{align*}
& \hspace{-0.5cm} X_{\nu}\defeq (X_{i_1}, \ldots, X_{i_k}),\\
&\hspace{-0.5cm} F^{\nu}(t)\defeq   \E \left( \I_{(\infty,t]} (Z) \big \vert X_{\nu} \right).
\end{align*}

\noindent In \cite{gamboa2015} the following Cram\'er von Mises sensitivity index is considered. Assume that $F$ is absolutely continuous. This normalized index (denoted  by  $C^{\nu}_{2}$) is defined as 

\begin{equation}
C^{\nu}_{2} \defeq \frac{N^{\nu}_{2,}}{\int_{\mathbb R} F(x)\left(1-F(x) \right) dF(x)} = 6 N^{\nu}_{2},
\end{equation}

where
\begin{align*}
N^{\nu}_{2} =\int  \E \left( \left[  F(t) -F^{\nu}(t) \right ]^2 \right) dF(t) \\
&\hspace{-4.7cm}= \int \left ( \E  \left[   \int \I_{(\infty,t]}(z)  d(F -F^{\nu})(z) \right ]^2\right )  dF(t)\\
&\hspace{-4.7cm}= \E \left\{ \E  \left(  \left[ \E \left(\I_{(\infty,t]}(Z) \right) - \E \left(\I_{(\infty,t]}(Z) \big \vert X_{\nu} \right)  \right]^2 \right) \right\}.
\end{align*}

\noindent To build  of our new index, we replace the  function $\I_{(\infty,t]}$  by the indicator function of the interval $[\min(s,t), \max(s,t)]$ (denoted by $h_{s,t}$),

\begin{displaymath}
h_{s,t}(x) \defeq \I_{ \{s \leq x \leq t\} }   +   \I_{\{t \leq x \leq s\}}=   \I_{ \{\min\{s,t\} \leq x \leq \max\{s,t\}\}}.  
\end{displaymath}

\noindent The function $h_{s,t}$, can be thought of as the indicator function of the ball of diameter $\overline{st}$ in $\mathbb{R}$. Further, let

\begin{equation}
\label{bla}
H(s,t)\defeq \E\left[ h_{s,t} (Z)\right] \quad \textrm{and}  \quad  H^{\nu}(s,t)= \E\left[ h_{s,t} (Z) \Big \vert X_{\nu} \right]. 
\end{equation}
\noindent  Obviously, $\E_{X_{\nu}}\left[ H^{\nu}(s,t) \right]= H(s,t)$. 

\begin{definition}
\label{sensibilidad}
The normalized ball sensitivity index (denoted by $B^{\nu}_{2}$) is then defined as,
  \begin{equation}
  B^{\nu}_{2} \defeq \frac{S^{\nu}_{2}}{\int_{\mathbb R^2} H(x,y)\left(1-H(x,y) \right) dF(y)dF(x)} = 6 S^{\nu}_{2},
  \end{equation}
  where 
  \begin{equation}
  S^{\nu}_{2} \defeq \E_{Z_1,Z_2} \left[ \E_{X_{\nu}} \left \{  \left[ H(Z_1,Z_2) - H^{\nu}(Z_1,Z_2)\right]^2 \right \} \right],
  \end{equation}
\noindent  and $Z_1$ and $Z_2$ are two independent copies of $Z$.
(The constant 6 follows from the well known fact that $\int_{t>s} (t-s)(1-(t-s))dtds = \int_{t\leq s} (t-s)(1-(t-s))dtds=1/12$.)

  \end{definition}

We may can rewrite $S^{\nu}_{2}$ as

\begin{displaymath}
S^{\nu}_{2}=\int_{ \M \times \M}  \E_{X_{\nu}} \left \{  \left( H(z_1, z_{2})-  H^{\nu}(z_1,z_{2})\right)^2 \right \} dF(z_1) dF(z_{2}).
\end{displaymath}

%

\subsection{Generalization for a  Riemannian manifold}

As discussed before, $\M$ is a Riemannian manifold of dimension $k$. So that, under some regularity conditions, given two points  $\{z_1, z_2\} \subset \M$, the following  function $h_{z_1, z_2 }: \M \rightarrow \{0,1\}$ is generically well defined,

\begin{displaymath}
h_{ z_1, z_2  }(t) \defeq \I_{ B_{z_1z_2}} (t).
\end{displaymath}

\noindent Let $Z_1,Z_2$ be independent copies of $Z$. Our normalized  sensitivity index is built by means of balls as in Definition \ref{sensibilidad}, 

 \begin{equation}
 \label{indic}
  B^{\nu}_{2} \defeq \frac{S^{\nu}_{2}}{D^{\nu}_{2}},
  \end{equation}
where 
\begin{displaymath}
S^{\nu}_{2} \defeq  \E_{Z_1,Z_2} \left[ \mbox{Var}_{X_{\nu}} \left \{ \E_Z \left[ h_{Z_1,Z_2}(Z) \Big \vert X_{\nu} \right] \right\} \right],
\end{displaymath}


and

\begin{displaymath}
D^{\nu}_{2} \defeq \E \left[ H(Z_1,Z_2) \left(1-H(Z_1,Z_2) \right) \right].
\end{displaymath}

\noindent  Notice that the index defined in the previous subsection is a particular case with  $\M=\R$.

\begin{remark}
 If   $B^{\nu}_{2}=0$ we have that  $\E_Z \left[ h_{z_1,z_2}(Z) \Big \vert X_{\nu} \right]= \E_Z \left[ h_{z_1,z_2}(Z) \right] \, a.s.$ for all $(z_1,z_2) \in \Omega \subset \mathcal{M}\times\mathcal{M}$ with $P(\Omega)=1$. Therefore, under the assumptions of Property \ref{propiedad1},  the probability measures  $P_{Z \vert X_{\nu}}$ and $P_{Z}$ are the same. 
\end{remark}

\subsection{Estimation}

In the next subsection, an estimator of the index $B^{\nu}_{2}$ is proposed.   The variance of the conditional mean is estimated using the pick and freeze method (see \cite{sobol1993} and \cite{sobol2001}). Further, as the expected value  $\E_{Z_1,Z_2}(\cdot)$  is a symmetric function of $Z_1$ and $Z_2$ we will use an $U$--statistic. At the end  of the next subsection, the steps for the construction of the estimator will be detailed.

\subsubsection{Estimation by \textit{``Pick and Freeze''} method}
The method consists in writing the variance of the conditional mean as a covariance. For $\nu \in \{1, \ldots, d \}$, let $X^{\nu}$ be the random vector that coincides with $X$ on its $\nu$ components and is independently regenerated on the others components. That is, $X^{\nu}_{\nu}=  X_{\nu}$  and $X^{\nu}_{\bar{\nu}}=  X^{'}_{\bar{\nu}}$, where $X^{'}$ is an independent copy of $X$ and $\bar{\nu}$  is the complementary set of $\nu$ ($\bar{\nu}= \{1, \ldots, d \} \setminus \nu $). We set
\begin{displaymath}
Z^{\nu} \defeq f(X^{\nu}).
\end{displaymath}

\noindent For sake of simplicity assume first that  $\E(Z)=0$. Then,

\begin{equation}
\mbox{Var}\left(  \E(Z \big \vert X_{\nu}) \right)= \E \left(  \E^2(Z \big \vert X_{\nu}) \right).
\label{pick}
\end{equation}

\noindent Since $Z$ and $Z^{\nu}$ are conditionally independent of $X_{\nu}$ (see \cite{klein2013}), so that

 \begin{align*}
 \mbox{cov} (Z, Z^{\nu})=  \E(Z Z^{\nu}) - \E(Z)\E( Z^{\nu})=  \E(Z Z^{\nu}) - \E^2(Z)\\
 &\hspace{-6.95cm} = \E \left[ \E(Z Z^{\nu} \big \vert X_{\nu}) \right] = \E \left[ \E( Z \big \vert X_{\nu})  \E( Z^{\nu} \big \vert X_{\nu}) \right] \\
  &\hspace{-6.95cm} =  \E \left[  \E^2  (Z \vert X_{\nu}) \right].
 \end{align*}

\noindent  Following  \cite{klein2013}, estimating the covariance  by a Monte Carlo Method we obtain the following  Pick Freeze estimator of  $\mbox{Var} \left[  \E  (Z \vert X_{\nu}) \right]$,

\begin{equation}
T^{\nu} \defeq \frac{1}{N}\sum_{j=1}^N Z_j Z_j^{\nu} - \left( \frac{1}{2N} \sum_{j=1}^N (Z_j + Z_j^{\nu}) \right)^2.
\end{equation}

\noindent Here, $Z_j$ and $Z^{\nu}_j$,  are $N$ independent copies of $Z$  and  $Z^{\nu}$ respectively.

\vspace{0.2 cm}
\underline{\emph{First step}.} Using the previous idea, we build  a consistent estimator of the numerator $S^{\nu}_{2}$. 
Let $\mathcal{P}_{N,p}$ be the set of all possible subsets of $\{1, \ldots,N\}$ with $p$ elements, that is,

\begin{displaymath}
\mathcal{P}_{N,p}= \{(i_1, \ldots, i_p) \in \{1, \ldots, N \}^p / i_1< \ldots < i_p  \}.
\end{displaymath}

\noindent Let $\tau=(i_1, \ldots, i_p) \in \mathcal{P}_{N,p}$. We denote by $\textbf{W}_{\tau}= \left( W_{i_1},\ldots ,W_{i_{p}} \right)$ a sample of $p$ independent copy of $Z$.
The estimator is obtained as follows,

\noindent  Further,  $D^{\nu}_{2}$  is  estimated analogously setting

\begin{align*}
\hat{D}^{\nu}_{2} \defeq \frac{1}{\binom{N}{2}} \sum_{\tau \in \mathcal{P}_{N,2}} \Big \{ \frac{1}{2N} \sum_{j=1}^N   \left( h_{\textbf{W}_{\tau}} \left( Z_j\right)+ h_{\textbf{W}_{\tau}} \left( Z_j^{\nu}\right) \right)  -  \\
& \hspace{-7.5cm}  - \left( \frac{1}{2N} \sum_{i=1}^N \left[ h_{\textbf{W}_{\tau}} \left( Z_i\right)+ h_{\textbf{W}_{\tau}} \left( Z_i^{\nu}\right) \right] \right)^2   \Big \}.
\end{align*}

\noindent The computation  of the  estimator is simple. Indeed, it only involves computation of indicator functions. Nevertheless, the computational time to perform  this estimator could be high. Indeed,  the number of sums is of  order $N^3$. Each term involves the determination of a geodesic. In the following section, we give some asymptotic properties of this estimator. We provide an exponential inequality that leads, from the Borel--Cantelli Lemma, to strong consistency.

\subsection{Asymptotic properties of  $\hat{B}^{\nu}_{2}$.}
We analyzed separately the strong convergence of $\hat{S}^{\nu}_{2}$ and $\hat{D}^{\nu}_{2}$. The strong convergence of  $\hat{S}^{\nu}_{2}$  follows from the next Theorem. Notice that another proof may be built using the fourth moment and the Rosenthal inequality for $U$--statistic of order $2$ developed in \cite{ibragimov1999}.

\textbf{Theorem (Exponential inequality )}
\label{exponencial}
Let $s>0$, there exist $N_0$ such that if $N>N_0$,
\begin{equation}
P \left( \left \vert \hat{S}^{\nu}_{2} -  S^{\nu}_{2} \right \vert > s   \right) \leq 16 \exp \left\{  -\frac{N \left(\frac{s}{9}\right)^2}{8} \right\}.
 \end{equation}

From the previous Theorem and the Borel--Cantelli Lemma the strong consistency of $\hat{S}^{\nu}_{2}$ follows.

\textbf{Corollary (Consistency of the estimator)}
\label{consistencia}
$\hat{S}^{\nu}_{2}$ is a consistent estimator of $S^{\nu}_{2}$.

Notice that we may analogously show the convergence of the denominator involved in (\ref{indic}).
\section{Simulations}
\label{simulaciones}
This simulation   section is made up on three examples. The first one shows the accuracy of the estimator  when the output is real valued. In the second example the output lies on a circle immersed in $\mathbb{R}^2$. The estimates are compared with those obtained in \cite{gamboa2015}. The last example shows that the index proposed in \cite{gamboa2015} may fail as sensitivity indicator.
 
\subsection{Example 1: Output on the real line}
\label{ejemplo1}
We study here an example where the Sobol index does not provide information on sensitivity, while $D^{\nu}_{2}$  and $B^{\nu}_{2}$ do. Let $Z=\alpha X_1 + X_2$  where $\alpha>0$, and $X_1,X_2$ are independent. Assume further that $X_1 \sim \textrm{Bernoulli}(p)$, $X_2 \sim F$, with $m=\E(X_2)$, $\sigma^2=\textrm{Var}(X_2)=\alpha^2 p(1-p)$.

We calculate $B^{1}_{2}$ in the Appendix while  $C^{1}_{2}$ has been calculated in \cite{gamboa2015}. In the case when $X_2 \sim U(0,b)$ with $b= \sqrt{12\alpha^2 p(1-p)}$. We find
\begin{displaymath}
B^{1}_{2}= 12p(1-p)\begin{cases}
 \left( \frac{\alpha}{b} \right)^3  \left(\frac{1}{3}- \frac{1}{4}\frac{\alpha}{b} \right) & \text{if } \alpha \leq b,\\
 1/12 & \text{if } \alpha > b,
\end{cases}
\end{displaymath}
\begin{displaymath}
C^{1}_{2}= 6p(1-p)\begin{cases}
 \left( \frac{\alpha}{b} \right)^2  \left(1- \frac{2}{3}\frac{\alpha}{b} \right) & \text{if } \alpha \leq b,\\
 1/3 & \text{if } \alpha > b,
\end{cases}
\end{displaymath}

\begin{figure}[ht]\centering 
\includegraphics[width=.3\textwidth]{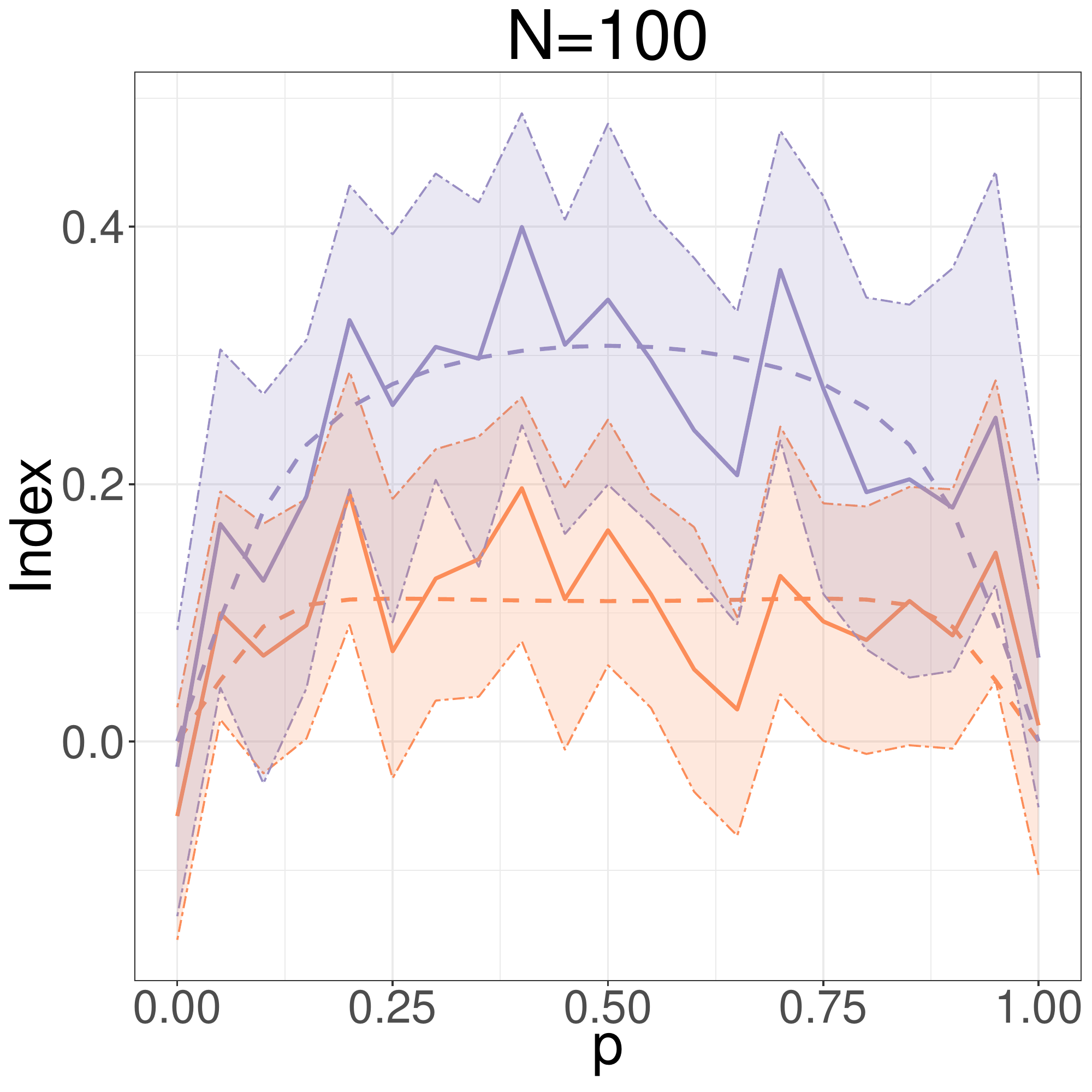}
\includegraphics[width=.3\textwidth]{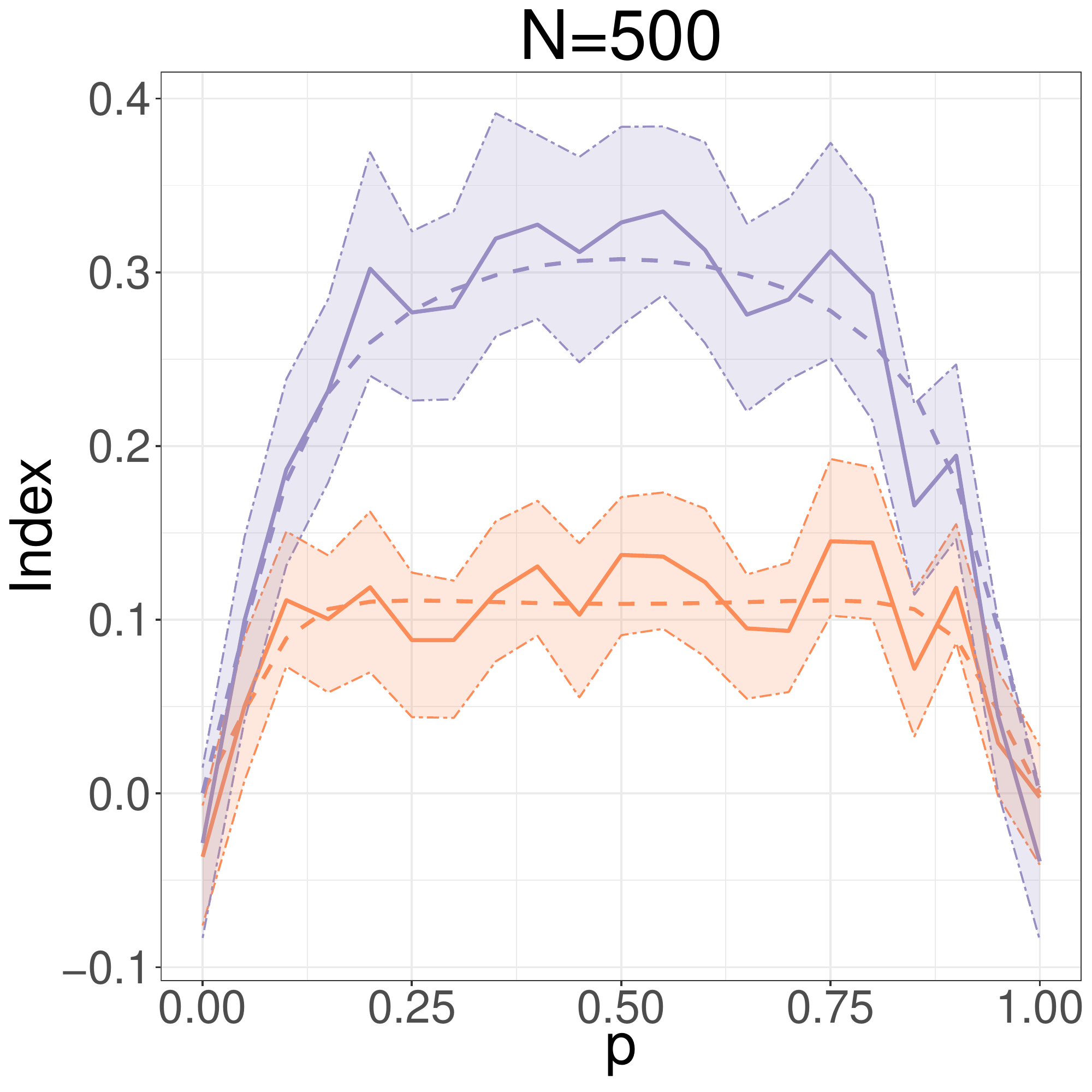}
\includegraphics[width=.3\textwidth]{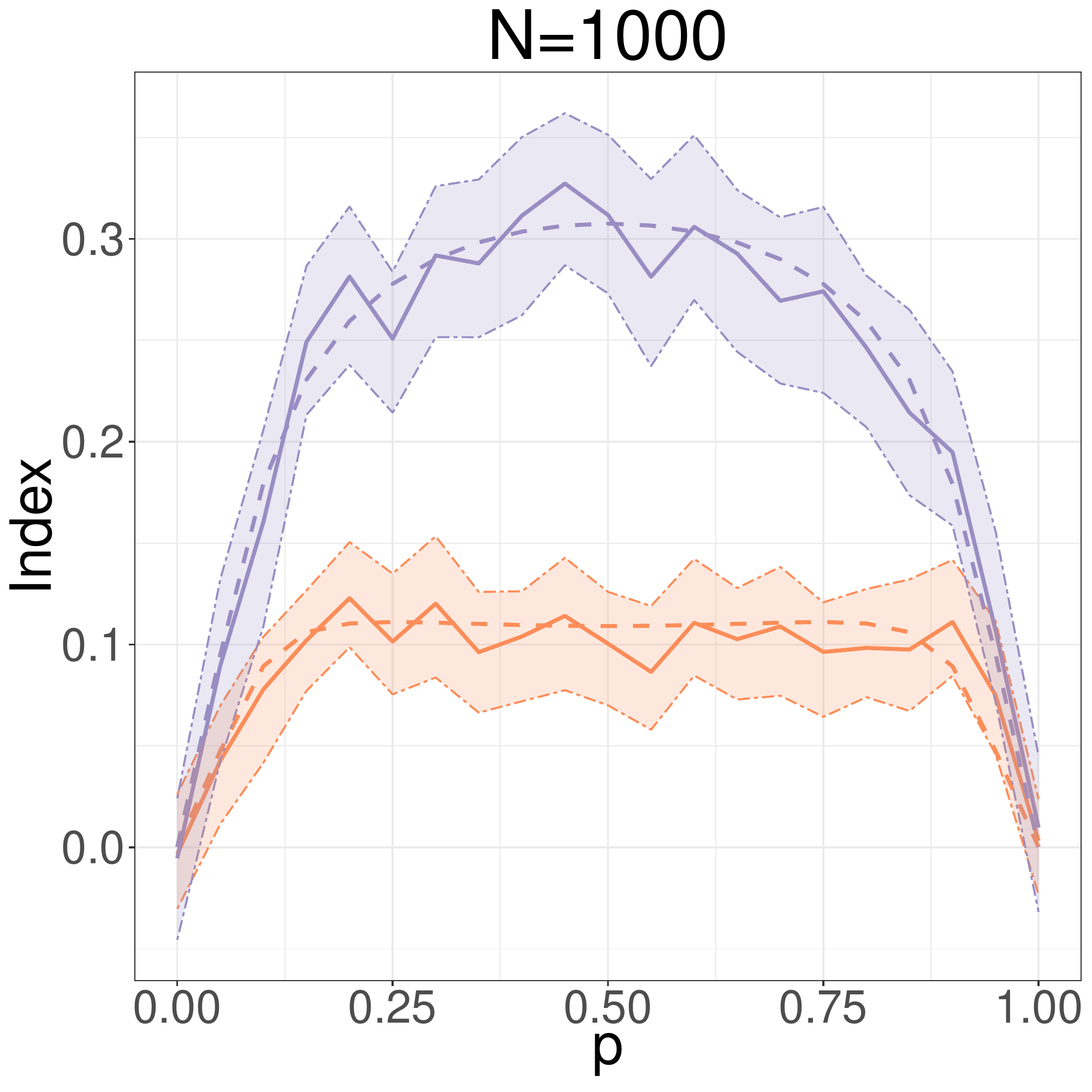}
\caption{Cram\'er von Mises indexes and our new indexes for Example \ref{ejemplo1}. In all cases,   ($\textendash \, \, \textendash$) are the true index  and ($\edge$) their estimate. Violet color depicts CVM sensitivity index while orange color depicts Ball sensitivity index. The $95\%$ bootstrap confidence intervals are represented with shading.  Left--hand Panel: $N=100$. Center--hand Panel: $N=500$. Right--hand Panel: $N=1000$.}
\label{index}
\end{figure}

In the Figure \ref{index}, the estimates of both indices are compared for different values of $p \in [0,1]$. Samples sizes $100$, $500$ and $1000$ are considered.  For each $p$, through resampling Bootstrap for $U$--statistics (see \cite{gine1992b}), $95\%$ confidence intervals are also given. In all cases, the mean squared deviation (MSD) of the estimators are computed.
The MSD of $\hat{B}^{1}_{2}$ is slightly smaller than the one for $\hat{C}^{1}_{2}$ (see Table \ref{MSDrecta}). Similar results are obtained when comparing  $\hat{B}^{2}_{2}$  and  $\hat{C}^{2}_{2}$.

\begin{table}[ht]
\centering
\begin{tabular}[t]{ccc}
\toprule
 Size &   $MSD_{\hat{B}^{1}_{2}}$ &  $MSD_{\hat{C}^{1}_{2}}$ \\
\midrule
$N=100$    & 0.051  & 0.067 \\
$N=500$    & 0.022 & 0.028 \\
$N=1000$   & 0.013  & 0.018  \\
\bottomrule
\end{tabular}
\caption{The MSD of the  $\hat{B}^{1}_{2}$ and $\hat{C}^{1}_{2}$ estimators for sample sizes $N=100,500$ and $1000$.}
\label{MSDrecta}
\end{table}%

We may conclude that the new proposed index has a similar behavior as the Cram\'er von Mises one for  the real line.

\subsection{Example 2: Output on a simple manifold immersed in $\mathbb{R}^2$} 
\label{ejemplo2}
The case where $Z$ ranges on the unit circle $S_1$ of  $\mathbb{R}^2$ is considered. We assume that the input vector has distribution  $$X=\binom{X_1}{X_2}\sim N \left[ \binom{\mu_1}{\mu_2},  \left(
\begin{array}{cc}
\sigma^2_1 & 0 \\
0 & \sigma^2_2
\end{array} \right)  \right].$$

Here, we study the case when the output $Z$ is the normalized version of $X$

\begin{displaymath}
\label{circu}
Z \defeq \frac{X}{\Vert X \Vert}.
\end{displaymath}

\begin{figure}[ht]\centering 
\includegraphics[width=.49\textwidth]{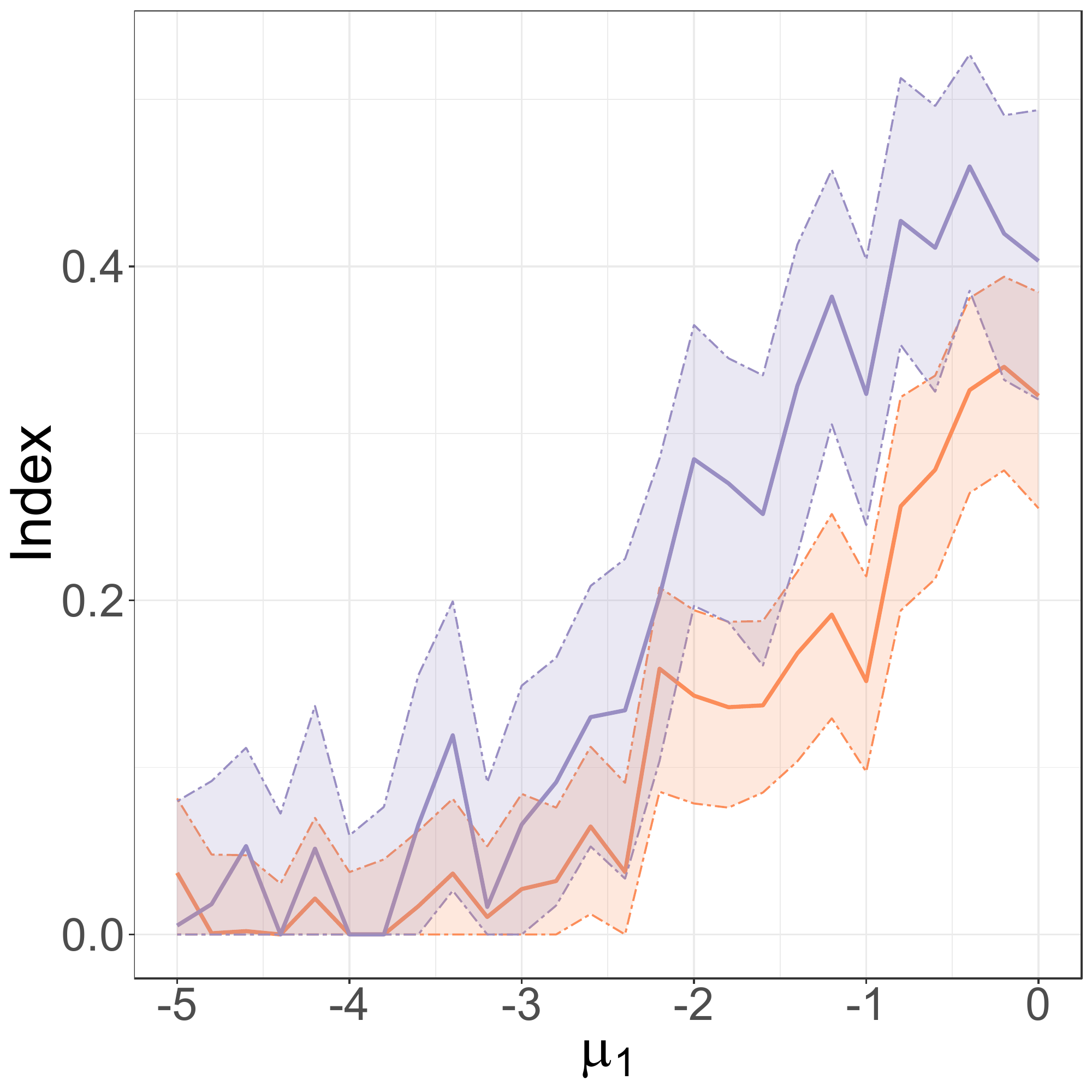}
\includegraphics[width=.49\textwidth]{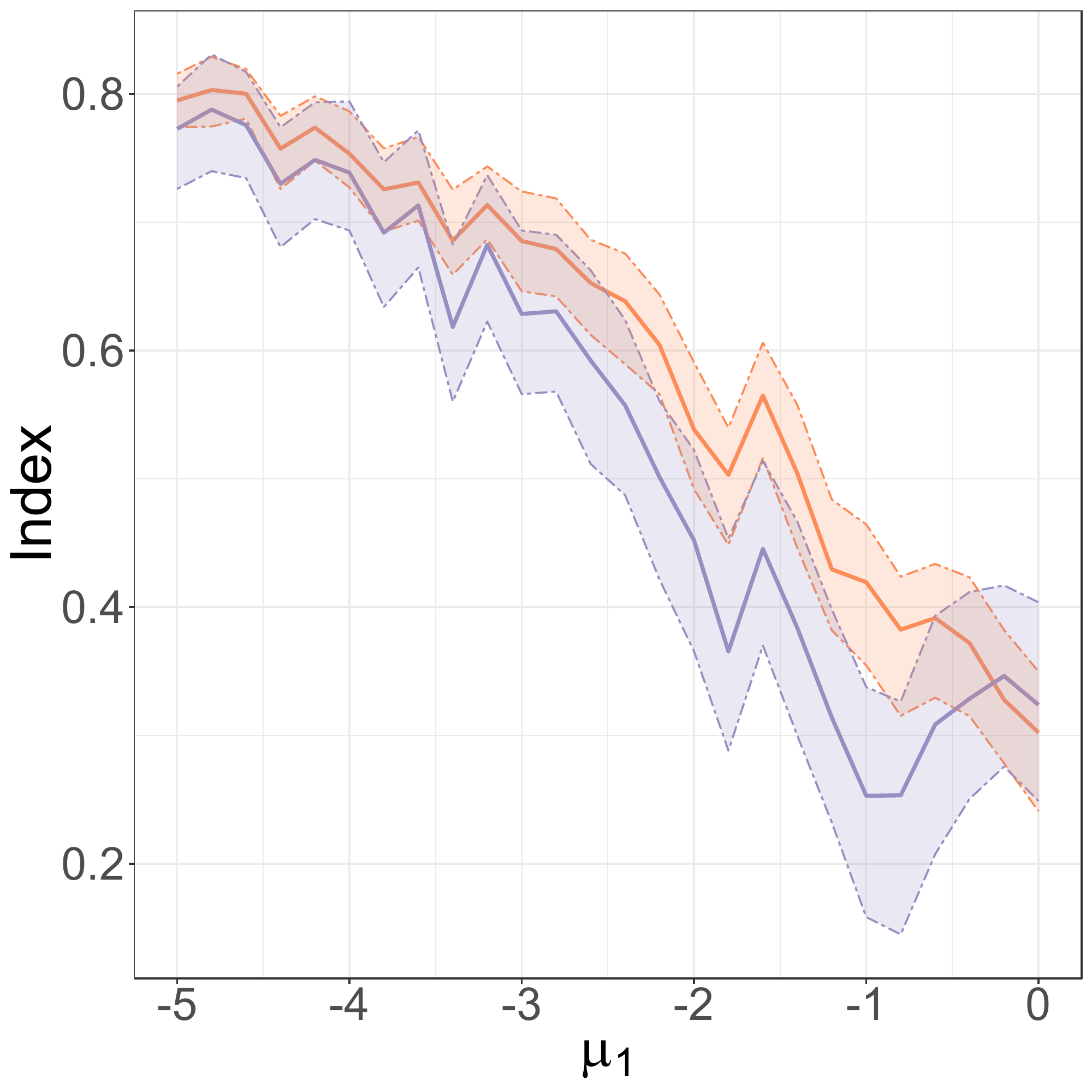}
\caption{Calculation of the indices $\hat{B}^{\nu}_{2}$ and $D^{\nu}_{2}$ in the equation \ref{circu} with $\mu_1 \in [-5,0]$ and $N=300$. For all cases, it is represented using  ($\edge$) the index estimation functions for $N=300$. In violet color depicts CVM sensitivity index and orange color depicts Ball sensitivity index. The $95\%$ bootstrap confidence intervals are represented with shading. Left--hand Panel: $\nu=1$. Right--hand Panel: $\nu=2$.}
\label{circulo}
\end{figure}

The distribution of $Z$ has been widely studied, we refer for example to \cite{mardia1972}, \cite{kendall1974} and \cite{watson1983}. In Figure \ref{circulo} the  index proposed  in \cite{gamboa2015} and the one studied here are depicted for the last system. These sensitivity indices are computed for $\mu_1 \in [-5,0]$. The other values are setted to  $\mu_2=0$ and $\sigma^2_1=\sigma^2_2=1$. We observe that the variability of $\hat{B}^{\nu}_{2}$ is  smaller that  the one of $\hat{D}^{\nu}_{2}$  for a sample size $N=300$.

\subsection{Example 3: Output on a manifold immersed in $\mathbb{R}^3$ }

Let us now consider a system with output in the following manifold $\M$ immersed in $\mathbb{R}^3$

\begin{displaymath}
\M= \left \{(x,y,z) \in \mathbb{R}^3 /  xyz=1, \, \, x,y,z >0  \right \}.
\end{displaymath}

\begin{figure}[ht]\centering 
\includegraphics[width=.80\textwidth]{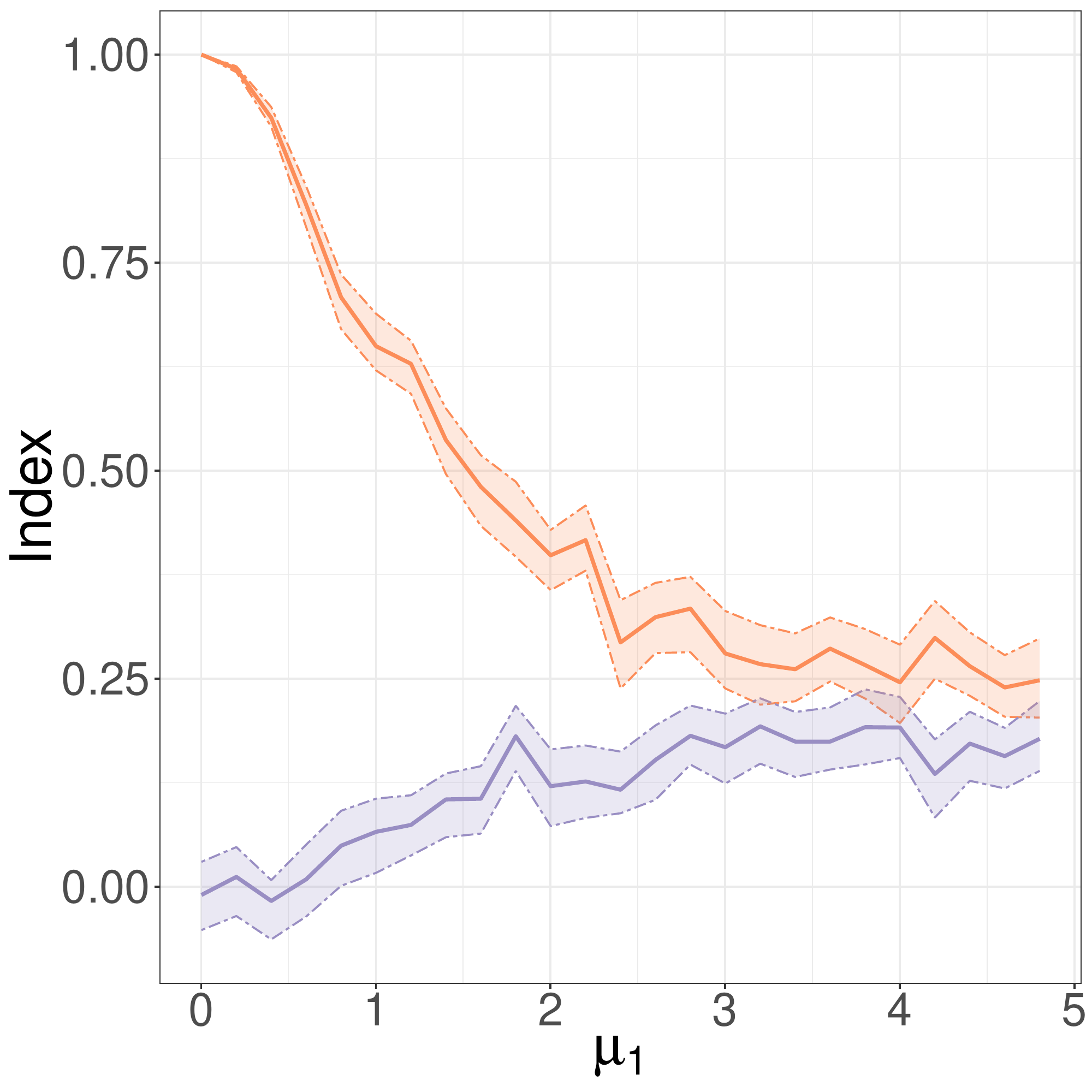}
\caption{Calculation of the index $\hat{B}^{\nu}_{2}$ }
\label{ej3}
\end{figure}

More precisely, let $Z=f(X,Y)\defeq (X+Y, \frac{1}{X}, \frac{X}{X+Y})$  where $X$ and $Y$ are iid random variables with Gamma distribution with parameters $(\mu_1,1)$, $(\mu_1>0)$. Since the function    $\I_{\{z \leq w\}} = 0$ for all $z, w \in M, z\neq w$, the index estimation $\hat{C}^{\nu}_{2}$ based on the left lower quadrants introduced in \cite{gamboa2015}, does not provide any information about the sensitivity. Also the second moment of $1/X$ does not exist, therefore it is not possible to compute the Sobol' index defined in \cite{klein2014}. By varying the parameter $\mu_1$, we calculate the estimation  of the ball sensitivity index $B^{1}_{2}$ and $B^{2}_{2}$ . For this purpose,  $N=1000$ pick-freeze samples were generated for each value of $Z$, $(Z_j, Z^{\nu}_j)$, $j=1,\ldots,N$; and other $1000$ samples of $Z$,    independent of $(Z_j, Z^{\nu}_j)$ for the corresponding $W_k$.
In Figure \ref{ej3} the values of the indices are observed by varying the parameter $\mu_1$ between  $0$ and $5$. It is clear from Figure \ref{ej3} that the new index is able to detect the effect of each variable on the output for different values of the parameter $\mu$.

\section{Sensitivity for isotropic matrix}
\label{datosreales}
We consider now the stiffness matrix $Z$ of an  isotropic materials as a function of the constants of Lam\'e  $\lambda$ and $\mu$ (see \cite{landau1965}, p. 13). Hence, if the temperature is constant we have 

\begin{equation*}
Z = \left(
\begin{array}{cccccc}
K+ 4\mu/3   & K-2\mu/3  & K-2\mu/3 & 0 & 0 & 0
 \\
K-2\mu/3  & K+ 4\mu/3  & K-2\mu/3 & 0 & 0 & 0 \\
K-2\mu/3  & K-2\mu/3  & K+ 4\mu/3 & 0 & 0 & 0 \\
0&0&0&\mu&0&0 \\
0&0&0&0&\mu&0 \\
0&0&0&0&0& \mu\\
\end{array} \right)
\end{equation*}

\noindent where $K= \lambda +2\mu/3$ is the volumetric modulus. The parameter $\mu$ is called the stiffness modulus, $K$ and $\mu$ take nonnegative values and are are modeled with a distribution supported by $\mathbb{R}^{+}$. The set of stiffness matrices is considered as a sub--manifold of the Riemannian manifold of the symmetric positive-definite matrices with the metric $g$. This manifold is usually denoted by $(\mathbb{P}_d,g)$, we refer to \cite{moakher2005} for more on the subject. Given two matrices $A$ and $B$  there exists an unique geodesic joining $A$ and $B$ given by,

\begin{equation}
\gamma(t) \defeq  A ^{1/2} \left(A ^{-1/2} B A ^{-1/2}  \right) ^ t A ^{1/2}. 
\end{equation}  

\noindent So that, we can calculate the midpoint between $A$ and $B$. We will denote by $A\#B$ this midpoint. We have

\begin{align}
A\#B = A ^{1/2} \left(A ^{-1/2} B A ^{-1/2}  \right) ^{1/2} A ^{1/2} \\
& \hspace{-6.5 cm} d(A,B)= \Vert \textrm{ln}  \left(A ^{-1/2} B A ^{-1/2}  \right) \Vert,
\end{align}
where $\Vert \cdot \Vert$  is the Hilbert--Schmidt norm. Note that the midpoint is simply the geometric mean of the matrices. 

\vspace{2 mm}

We will focus on the case when  $K$ and $\mu$ are independent random variables and on two scenarios:  

\begin{tabular}{lll}
\textbf{Case 1} & $\quad K \sim \gamma( 1/\lambda_K, \lambda_K )$ & $ \mu \sim \gamma(1/\lambda_{\mu}, \lambda_{\mu} )$ \\ 
\textbf{Case 2} & $\quad K \sim \textrm{U}(1-\lambda_K , 1+\lambda_K )$ & $\mu \sim \textrm{U}(1-\lambda_{\mu},1+\lambda_{\mu})$ \\ 
\end{tabular} 
 
Tables \ref{indice_ug1} and \ref{indice_ug2}  show the values of  $\hat{B}^{1}_{2}$ and $\hat{B}^{2}_{2}$  for  various values of $(\lambda_K,\lambda_{\mu})$. The sample size was $N=500$. In both tables we can see that the new index provides significant information regarding the sensibility of $Z$ with respect to the input variables $K$ and $\mu$.  In particular, when the variance of the input variables increases (larger $\lambda_{\mu}$ or $\lambda_{K}$) the value of the index increases, i.e. a a greater effect is observed on the output.

\begin{table}[ht]
\begin{center}
\caption{Index on isotropic matrix with Gamma distributions parameter}
\scalebox{0.8}{
\begin{tabular}{c | cccc|| cccc }
  \toprule[0.6 mm]
  Distribution  &  \multicolumn{4}{c}{Case 1: $B^1$ } & \multicolumn{4}{c}{Case 1: $B^2$ }             \\
\toprule[0.3 mm]
 $ \lambda_{\mu} \Big \backslash \lambda_K$ & \textbf{0.001} & \textbf{0.01}  & \textbf{0.1} & \textbf{1} &\textbf{0.001} & \textbf{0.01}  & \textbf{0.1} & \textbf{1} \\
  \toprule[0.3 mm]
  \textbf{0.001} & 0.625   & 0.212 & 0.016 & 0.001           & 0.083 & 0.435 & 0.855 & 0.980  \\
  \textbf{0.01} &0.925 & 0.593 & 0.215 & 0.033           & 0.004 & 0.072 & 0.458 & 0.865  \\
 \textbf{0.1} & 0.987 & 0.912 & 0.587 & 0.184          & 0.001 &0.006& 0.137 & 0.518         \\
\textbf{1} & 0.999 & 0.990 & 0.930 & 0.600       & 0.000 & 0.007 & 0.210 & 0.311   \\

  \toprule[0.6 mm]
\end{tabular}}
\label{indice_ug1}
\end{center}
\end{table}

\begin{table}[ht]
\begin{center}
\caption{ Index on isotropic matrix with Uniform distributions parameter}
\scalebox{0.8}{
\begin{tabular}{c | cccc || cccc }
  \toprule[0.6 mm]
  Distribution  &  \multicolumn{4}{c}{Case 1: $B^1$ } & \multicolumn{4}{c}{Case 1: $B^2$ }             \\
\toprule[0.3 mm]
  $ \lambda_{\mu} \Big \backslash \lambda_K$ & \textbf{0.001} & \textbf{0.01}  & \textbf{0.1} & \textbf{1} &\textbf{0.001} & \textbf{0.01}  & \textbf{0.1} & \textbf{1} \\
  \toprule[0.3 mm]
  \textbf{0.001} & 0.620   & 0.008  & 0.001 & 0.000          & 0.109 & 0.848 & 0.997 & 1.000  \\
  \textbf{0.01}  &0.989     & 0.623 & 0.003 & 0.001           & 0.018 & 0.092 & 0.849 & 0.997  \\
 \textbf{0.1} & 1.000      & 0.989 & 0.623 & 0.624          & 0.016 &0.016& 0.102 & 0.846        \\
\textbf{1} & 1.000         & 1.000 & 0.990 & 0.987            & 0.015 & 0.016 & 0.100 & 0.211     \\

  \toprule[0.6 mm]
\end{tabular}}
\label{indice_ug2}
\end{center}
\end{table}

\section{Conclusions and Future Work}
\label{conclusiones}

We introduce a new framework to measure the  global sensitivity  for system outputs valued on a Riemannian manifolds. This sensitivity index has nice properties:
\begin{enumerate}
\item It is a generalization of the one proposed by \cite{gamboa2015}.
\item  It is built on the whole distributions and not only on moments.

\item The construction of the index lies on the geometry of the support of the output through the  geodesic distance.

\item A pick-freeze like estimator of the index based on an $U$--statistic is  proposed. It is easy to compute. The consistency of the estimator is shown. 

\item Considering  various simulation scenarios, the desirable properties of the estimator are illustrated. Using a Bootstrap resampling method for $U$--statistics we construct confidence intervals.

\item  The sensitivity index for an isotropic matrix model is studied.

\item In future research we will challenge to study new applications, for instance, the impact of inputs generated by a dynamic system, when the output is supported by a Riemannian Lie  Group.

\end{enumerate}

%

\section*{Appendix}
\begin{enumerate}[label=(\Alph*)]
\item \textbf{Proof of Property \ref{propiedad1}.}

In \cite{christensen1970} it is shown that both the sets

\begin{equation}\label{determ}
 \mathcal B \defeq \{{B}_{p,r} \big/ p \in \M  \, \textrm{and} \, r>0 \}, \ \mbox{and} \ \mathcal B_{\delta} \defeq \{{B}_{p,r} \big/ p \in \M  \, \textrm{and} \, 0<r<\delta \},
 \end{equation}
are determining classes. We aim to  show that the family of balls $$\mathcal{B}_{pq}:= \{B_{pq}, p,q \in \M\ \textrm{and there is a unique geodesic between $p$ and $q$} \},$$ is also a determining class if the radius of injectivity  $r_{iny}$ of $\M$ is positive.
Recall that in a metric space the family of closed sets $\mathcal{C}$ is a determining class (\cite{billingsley1999}, page $7$). Therefore, it suffices to show that if two probabilities $\nu$ and $\eta$ coincides on $\mathcal{B}_{pq}$, then $\nu(A)=\eta(A)$ for all $A \in \mathcal{C}$. Let  $\partial A$ the topological boundary of $A$ and $\partial A^{\epsilon} \defeq \bigcup_{x \in \partial A} B(x, \epsilon)$ for $\epsilon >0$. So,

$$\nu (A)=  \nu \left( A \setminus (\partial A^{\epsilon} \cup \partial \M^{\epsilon} ) \right) + \nu \left( A \cap (\partial A^{\epsilon} \cup \partial \M^{\epsilon} ) \right).$$

Let  $\epsilon^{*} \defeq \min (\epsilon/2, r_{iny})$.  Since by (\ref{determ}) $\mathcal B_{\epsilon^{*}}$ is a determining class, using the exponential map  we derive that the set of balls in $\mathcal B_{\epsilon^{*}}$ 
that determine the $\nu$-probability of $A \setminus (\partial A^{\epsilon^{*}} \cup \partial \M^{\epsilon^{*}})$  are in
  $\mathcal{B}_{pq}$. Therefore, $$\nu \left( A \setminus (\partial A^{\epsilon^{*}} \cup \partial \M^{\epsilon^{*}} ) \right)= \eta \left( A \setminus (\partial A^{\epsilon^{*}} \cup \partial \M^{\epsilon^{*}}) \right).$$ 
Finally, from the dominated convergence theorem  and the assumptions, we conclude that
$$\nu (A)=  \eta(A)- \eta(\partial A  \cup \partial \M) + \nu(\partial A  \cup \partial \M) =\eta(A).$$


\QEDB

\item \textbf{Proof of Property \ref{exponencial}.}


We will show that given $s>0$, there is an $N_0$ such that for any $N>N_0$,
\begin{equation}
\label{exponencial2}
P \left( \left \vert \hat{S}^{\nu}_{2} -  S^{\nu}_{2} \right \vert > 9s   \right) \leq 16 \exp \left\{  -\frac{N^2}{8} \right\}.
 \end{equation}
 
We will prove that $$P \left(  \hat{S}^{\nu}_{2} -  S^{\nu}_{2}  > 9s   \right)  \leq 8 \exp \left\{  -\frac{Ns^2}{8} \right\},$$  and it is analogous for the other tail.

For $1 \leq j,k \leq N$ and $\tau \in \mathcal{P}_{N,2}$ let, 
\begin{itemize}
\item $\textbf{W}_{\tau}= \left( W_{k_1},W_{k_2} \right)$  
\item $\textbf{Z}_{j}= \left( Z_j, Z_j^{\nu}  \right)$
\item $G(\textbf{Z}_{j},\textbf{W}_{\tau})= h_{\textbf{W}_{\tau}} (Z_j) h_{\textbf{W}_{\tau}} (Z_j^{\nu})$
\item $J(\textbf{Z}_{j},\textbf{W}_{\tau})= \frac{1}{2} \left[  h_{\textbf{W}_{\tau}} (Z_j)+ h_{\textbf{W}_{\tau}} (Z_j^{\nu}) \right]$
\item  $H(\textbf{Z}_{i},\textbf{Z}_{j},\textbf{W}_{\tau})=J(\textbf{Z}_{i},\textbf{W}_{\tau})J(\textbf{Z}_{j},\textbf{W}_{\tau})$
\end{itemize}

The proof is built on three steps:

\begin{itemize}
\item \textbf{Step 1} [Rewrite the difference $\hat{S}^{\nu}_{2}-S^{\nu}_{2}$]
Similarly to  \cite{gamboa2015}  but for an $U$--statistics, we have that

\begin{align*}
\hat{S}^{\nu}_{2} = \frac{1}{N\binom{N}{2}} \sum_{\substack{j \in \{1,\ldots, N\} \\ \tau \in \mathcal{P}_{N,2} }} G(\textbf{Z}_{j},\textbf{W}_{\tau}) -\frac{1}{N^2\binom{N}{2}} \sum_{\substack{\{i,j\} \in \{1,\ldots, N\} \\ \tau \in \mathcal{P}_{N,2} }} H(\textbf{Z}_{i},\textbf{Z}_{j},\textbf{W}_{\tau})  \\
&\hspace{-11.7cm} =\frac{1}{N\binom{N}{2}} \sum_{\substack{j \in \{1,\ldots, N\} \\ \tau \in \mathcal{P}_{N,2} }} \left \{ G(\textbf{Z}_{j},\textbf{W}_{\tau}) -\E \left[ G(\textbf{Z}_{j},\textbf{W}_{\tau}) \right] \right\}-\\
 &\hspace{-11.7cm}  -\frac{1}{N^2\binom{N}{2}} \sum_{\substack{\{i,j\} \in \{1,\ldots, N\} \\ \tau \in \mathcal{P}_{N,2} }} \left \{ H(\textbf{Z}_{i},\textbf{Z}_{j},\textbf{W}_{\tau}) -\E\left[ H(\textbf{Z}_{i},\textbf{Z}_{j},\textbf{W}_{\tau}) \right] \right \}  +  \\
&\hspace{-11.7cm} + \E \left[ G(\textbf{Z}_{1},\textbf{W}_{\tau_1}) \right] - \left( 1-\frac{1}{N}\right) \E\left[ H(\textbf{Z}_{1},\textbf{Z}_{2},\textbf{W}_{\tau_1}) \right]  -   \frac{1}{N}\E\left[ H(\textbf{Z}_{1},\textbf{Z}_{1},\textbf{W}_{\tau_1}) \right],
\end{align*}
and
\begin{align*}
S^{\nu}_{2}= \E_{\textbf{W}_1} \left\{ \textrm{Var}_{X_{\nu}}\left( H^{\nu}(\textbf{W}_1)\right) \right \}&\\
&\hspace{-4.5cm}  =  \E_{\textbf{W}_1} \left\{  \textrm{Var}_{X_{\nu}}\left( \E_Z (h_{\textbf{W}_1} (Z)/ X_{\nu}) \right) \right \}=  \E_{\textbf{W}_{1}} \left\{  \textrm{cov} \left( h_{\textbf{W}_{1}} (Z_1), h_{\textbf{W}_{1}} (Z_1^{\nu}) \right) \right \} \\
&\hspace{-4.5cm}= \E_{\textbf{W}_{1}} \left\{  \E_{\textbf{Z}_1} \left( h_{\textbf{W}_{1}} (Z_1) h_{\textbf{W}_{1}} (Z_1^{\nu}) \right) \right \}-  \E_{\textbf{W}_{1}}\left\{ \left[  \E_{\textbf{Z}_1}  \left( h_{\textbf{W}_{1}} (Z_1) \right)\right]^2 \right \}.
\end{align*}

%
%

Now, we may  decompose the error   into three terms
\begin{align*}
\hat{S}^{\nu}_{2} -  S^{\nu}_{2}= \underbrace{\frac{1}{N\binom{N}{2}} \sum_{\substack{j \in \{1,\ldots, N\} \\ \tau \in \mathcal{P}_{N,2} }} \left \{ G(\textbf{Z}_{j},\textbf{W}_{\tau}) -\E \left[ G(\textbf{Z}_{j},\textbf{W}_{\tau}) \right] \right\}}_{(A)}-\\
 &\hspace{-9cm} \underbrace{ -\frac{1}{N^2\binom{N}{2}} \sum_{\substack{\{i,j\} \in \{1,\ldots, N\} \\ \tau \in \mathcal{P}_{N,2} }} \left \{ H(\textbf{Z}_{i},\textbf{Z}_{j},\textbf{W}_{\tau}) -\E\left[ H(\textbf{Z}_{i},\textbf{Z}_{j},\textbf{W}_{\tau}) \right] \right \}  }_{(B)} + \\
&\hspace{-9cm} \underbrace{+\frac{1}{N}  \left \{ \E\left[ H(\textbf{Z}_{1},\textbf{Z}_{2},\textbf{W}_{\tau_1}) \right]  - \E\left[ H(\textbf{Z}_{1},\textbf{Z}_{1},\textbf{W}_{\tau_1}) \right] \right \}.}_{(C)}
\end{align*}

In the following steps, we bounded the terms $(A)$ and $(B)$ for a sufficiently large $N$. 

In  steps  $2$ and $3$ we use Hoeffding inequality (see for example \cite{hoeffding1963}) for independent and real variables. That is, if  $W_1, \ldots, W_N$ are  independent centred random variables, supported by $[a_i,b_i] \quad i=1, \ldots,N$,
\begin{equation}
P\left( \sum_{i=1}^{N} W_i > s \right) \leq  \exp \left\{ - \frac{2s^2}{\sum_{i=1}^{N}(b_i-a_i)^2} \right\} \quad (s>0).
\end{equation}

We also use an extension of the previous inequality for $U$--statistics of order $2$, see \cite{serfling1980}, Theorem A, (5.6).  If  $s=s(X_1, X_2)$ is the $U$--statistics kernel of $U_n$ such that  $\E(s(X_1,X_2))=\theta$ and  $a \leq s(x_1, x_2) \leq b$. So, for  $s>0$ and  $N>2$,
\begin{equation}
P\left( U_n -\theta > s \right) \leq  \exp \left\{ - \frac{Ns^2}{(b-a)^2} \right\}
\end{equation}
For us, the kernels are centred bounded indicator functions.

\item \textbf{Step 2} [Bounds for $(A)$.]

Let $G_c(\textbf{Z}_{j},\textbf{W}_{\tau})= G(\textbf{Z}_{j},\textbf{W}_{\tau}) -\E \left[ G(\textbf{Z}_{j},\textbf{W}_{\tau}) \right]$ and $\tilde{G}(Z_i, \textbf{W}_{\tau})= G_c(Z_i, \textbf{W}_{\tau})  - \E_{Z_i} (G_c (Z_i, \textbf{W}_{\tau}))$,
\begin{displaymath}
A= \underbrace{\frac{1}{N\binom{N}{2}} \sum_{\substack{ i \in \{1,\ldots, N\} \\ \tau \in \mathcal{P}_{N,2} }}  \tilde{G}(Z_i,\textbf{W}_{\tau})}_{A_1} + \underbrace{\frac{1}{\binom{N}{2}} \sum_{\tau \in \mathcal{P}_{N,2} } \E_Z \left[ G_c(Z,\textbf{W}_{\tau}) \right] }_{A_2}
\end{displaymath}

\begin{align*}
P(A_1>s)= P \left(\frac{1}{N\binom{N}{2}} \sum_{\substack{ i \in \{1,\ldots, N\} \\ \tau \in \mathcal{P}_{N,2} }}  \tilde{G}(Z_i,\textbf{W}_{\tau}) >s \right) \\
&\hspace{-7 cm} = P \left( \sum_{i=1}^N  \frac{1}{\binom{N}{2}}\sum_{\tau \in \mathcal{P}_{N,2} } \tilde{G}(Z_i,\textbf{W}_{\tau}) >Ns\right)\\
&\hspace{-7 cm} \leq \exp \left\{  -\frac{s^2N}{2} \right\} \leq  \exp \left\{ -\frac{Ns^2}{8} \right\},
\end{align*}
while
\begin{align*}
P(A_2>s)= P \left( \frac{1}{\binom{N}{2}} \sum_{\tau \in \mathcal{P}_{N,2} } \E_Z \left[ G_c(Z,\textbf{W}_{\tau}) \right] >s \right) \\
&\hspace{-7 cm} \leq   \exp \left\{\frac{-2Ns^2}{8} \right\}\leq \exp \left\{  -\frac{Ns^2}{8} \right\}.
\end{align*}


But, $\left\{ A_1 + A_2> 2s \right\} \subset \{A_1>s\} \cup \{A_2>s \} $, so
\begin{equation}
P(A>2s) \leq 2\exp \left\{  -\frac{Ns^2}{8} \right\}.
\end{equation}

\item \textbf{Step 3} [Bounds for $(B)$.]

Let $H_c(\textbf{Z}_{i},\textbf{Z}_{j},\textbf{W}_{\tau}) = H(\textbf{Z}_{i},\textbf{Z}_{j},\textbf{W}_{\tau}) -\E\left[ H(\textbf{Z}_{i},\textbf{Z}_{j},\textbf{W}_{\tau}) \right]$ and  $\tilde{H}(Z_i, Z_j , \textbf{W}_{\tau})= H_c (Z_i,Z_j, \textbf{W}_{\tau})  - \E_{Z_j} (H_c (Z_i, Z_j,\textbf{W}_{\tau}) ) $, 

\begin{align*}
B= \underbrace{\frac{1}{N^2\binom{N}{2}} \sum_{\substack{\{i,j\} \in \{1,\ldots, N\} \\ \tau \in \mathcal{P}_{N,2} }}  \tilde{H}(Z_i,Z_j, \textbf{W}_{\tau})}_{B_1}+ \\
&\hspace{-6.5cm}  + \underbrace{\frac{1}{N\binom{N}{2}} \sum_{\tau \in \mathcal{P}_{N,2} } \E_Z \left[ H_c(Z,Z,\textbf{W}_{\tau}) \right] }_{B_2} + \\
&\hspace{-6.5 cm}+ \underbrace{\frac{N-1}{N^2\binom{N}{2}} \sum_{\substack{ j \in \{1,\ldots, N\} \\ \tau \in \mathcal{P}_{N,2} } } \E_Z \left[ H_c(Z,Z_j,\textbf{W}_{\tau}) \right] }_{B_3}.
\end{align*}

We have $\left\{ \sum_{i=1}^3B_i> 6s \right\} \subset \{B_1>3s\} \cup \{B_2>s \}\cup \{B_3>2s \} $.




So that, 
\begin{align*}
P(B_1>3s) \leq P \left(\sum_{j} \frac{1}{\binom{N}{2}}  \sum_{ \tau \in \mathcal{P}_{N,2} }  \tilde{H}(Z_j,Z_j,\textbf{W}_{\tau}) >N^2s \right)+ \\
&\hspace{-8.5 cm}+ P \left( \sum_{j>i} \frac{1}{\binom{N}{2}}  \sum_{ \tau \in \mathcal{P}_{N,2} }   \tilde{H}(Z_i,Z_j,\textbf{W}_{\tau}) >N^2s  \right)+ \\
&\hspace{-8.5 cm}+ P \left( \sum_{j<i} \frac{1}{\binom{N}{2}} \sum_{ \tau \in \mathcal{P}_{N,2} }   \tilde{H}(Z_i,Z_j,\textbf{W}_{\tau}) >N^2s \right)  \leq \\
&\hspace{-8.5cm} \leq \exp \left\{  -\frac{s^2N^3}{2} \right\} + 2 \exp \left\{  \frac{-Ns^2N^2}{4} \right\}\leq 3\exp \left\{  -\frac{Ns^2}{8} \right\}.
\end{align*}
Further,
\begin{align*}
P(B_2>s)  \leq P \left( \frac{1}{\binom{N}{2}} \sum_{\tau \in \mathcal{P}_{N,2} } \E_Z \left[ H_c(Z,Z,\textbf{W}_{\tau}) \right]  > Ns \right) \\
&\hspace{-8 cm} \leq  \exp \left\{  -\frac{Ns^2N^2}{4} \right\} \leq \exp \left\{  -\frac{Ns^2}{8} \right\}.
\end{align*}

Finally,
\begin{align*}
P(B_3>2s) \leq P \left( \frac{N-1}{N^2\binom{N}{2}} \sum_{\substack{ j \in \{1,\ldots, N\} \\ \tau \in \mathcal{P}_{N,2} } } \E_Z \left[ H_c(Z,Z_j,\textbf{W}_{\tau}) \right] >2s \right)  \\
&\hspace{-10.1 cm} \leq P \left( \sum_{ j \in \{1,\ldots, N\}  } \frac{1}{\binom{N}{2}} \sum_{\tau \in \mathcal{P}_{N,2} }\E_Z \left[ H_c(Z,Z_j,\textbf{W}_{\tau}) \right] >2s \frac{N^2}{N-1}\right) \\
&\hspace{-10.1 cm}  \leq P \Bigg( \sum_{ j \in \{1,\ldots, N\}  } \frac{1}{\binom{N}{2}} \sum_{\tau \in \mathcal{P}_{N,2} } \Big \{ \E_Z \left[ H_c(Z,Z_j,\textbf{W}_{\tau}) \right]- \\
&\hspace{-10.1 cm} -\E_{Z,Z_j} \left[ H_c(Z,Z_j,\textbf{W}_{\tau}) \right] \Big \} >\frac{s N^2}{N-1}\Bigg) + \\
&\hspace{-10.1 cm} + P \left(  \frac{1}{\binom{N}{2}} \sum_{\tau \in \mathcal{P}_{N,2} }\E_{Z_1,Z_2} \left[ H_c(Z_1,Z_2,\textbf{W}_{\tau}) \right]  >s \frac{N}{N-1}\right) \\
&\hspace{-10.1 cm} \leq  \exp \left\{  -\frac{s^2N^4}{2N(N-1)^2} \right\} + \exp \left\{  -\frac{2Ns^2 N^2}{8(N-1)^2} \right\} \leq 2\exp \left\{  -\frac{Ns^2}{8} \right\}.
\end{align*}

Collecting the previous inequality we obtain,
\begin{equation}
P(B>6s) \leq 6\exp \left\{  -\frac{Ns^2}{8} \right\}.
\end{equation}
\end{itemize}

\noindent  Now, using the bounds obtained in steps $2$ and $3$, we may conclude that,  for large enough  large $N$,  inequality (\ref{exponencial2}) holds.

\QEDB

\item \textbf{Proof to Example \ref{ejemplo1}.}

Let  $P(z_1, z_2)= F(z_2)-F(z_1)$, then we have 
 \begin{align*}
H(s,t)=\I_{t <s} +F_Z(t) -F_Z(s))\\
 &\hspace{-3.4cm} =  \I_{t <s} +(1-p)P(s,t)+pP(s-\alpha,t-\alpha),
 \end{align*}
 and,
\begin{displaymath}
H^1(s,t)=\begin{cases}
\I_{t <s} + P(s,t) & \text{if } X_1=0,\\
\I_{t <s} +P(s-\alpha,t-\alpha) & \text{if } X_1=1.
\end{cases}
\end{displaymath}
Therefore
\begin{small}
\begin{displaymath}
H(s,t)-H^1(s,t)=\begin{cases}
p \left[ P(s,t)-P(s-\alpha,t-\alpha)\right] & \text{if } X_1=0,\\
(1-p)\left[P(s-\alpha,t-\alpha)- P(s,t)\right]& \text{if } X_1=1,
\end{cases}
\end{displaymath}
\begin{displaymath}
\E_1 \left[\left(H(s,t)-H^1(s,t)\right)^2 \right] =p(1-p)\left[ P(s,t)-P(s-\alpha,t-\alpha)\right]^2,
\end{displaymath}
\end{small}
and,
\begin{align*}
S^{1}_{2}=p(1-p) \E\left[ \left(P(Z_1-\alpha,Z_1)-P(Z_2-\alpha,Z_2) \right)^2 \right]\\
 &\hspace{-6.9cm} = 2p(1-p) \textrm{Var} \left[ F(Z)-F(Z-\alpha) \right].
\end{align*}

Let  $X_2 \sim U(0,b)$, with $b= \sqrt{12\alpha^2 p(1-p)}$. We  obtain

\begin{displaymath}
\E\left(  F(Z)- F(Z-\alpha)  \right)=\begin{cases}\lim_{N\rightarrow \infty} \frac{1}{N^2}
  \frac{\alpha}{b}   \left(1- \frac{1}{2}\frac{\alpha}{b} \right) & \text{if } \alpha \leq b,\\
 1/2 & \text{if } \alpha > b.
\end{cases}
\end{displaymath}

\begin{displaymath}
\E\left[ \left( F(Z)- F(Z-\alpha) \right)^2 \right]=\begin{cases}
 \left( \frac{\alpha}{b} \right)^2  \left(1- \frac{2}{3}\frac{\alpha}{b} \right) & \text{if } \alpha \leq b,\\
 1/3 & \text{if } \alpha > b.
\end{cases}
\end{displaymath}

\begin{displaymath}
\textrm{Var} \left( F(Z)- F(Z-\alpha) \right)=\begin{cases}
 \left( \frac{\alpha}{b} \right)^3  \left(\frac{1}{3}- \frac{1}{4}\frac{\alpha}{b} \right) & \text{if } \alpha \leq b,\\
 1/12 & \text{if } \alpha > b.
\end{cases}
\end{displaymath}

So,

\begin{displaymath}
D^{1}_{2,CVM}= p(1-p)\begin{cases}
 \left( \frac{\alpha}{b} \right)^2  \left(1- \frac{2}{3}\frac{\alpha}{b} \right) & \text{if } \alpha \leq b,\\
 1/3 & \text{if } \alpha > b,
\end{cases}
\end{displaymath}

and

\begin{displaymath}
S^{1}_{2}= 2p(1-p)\begin{cases}
 \left( \frac{\alpha}{b} \right)^3  \left(\frac{1}{3}- \frac{1}{4}\frac{\alpha}{b} \right) & \text{if } \alpha \leq b,\\
 1/12 & \text{if } \alpha > b.
\end{cases}
\end{displaymath}

\end{enumerate}

\QEDB


\begin{thebibliography}{10}
\providecommand{\url}[1]{\normalfont{#1}}
\providecommand{\urlprefix}{Available from: }

\bibitem{tarantola2008}
Rocquigny~ED, Devictor~N, Tarantola~S. Uncertainty in industrial practice.
  Wiley Online Library; 2008.

\bibitem{saltelli2000}
Saltelli~A, Chan~K, Scott~E. Sensitivity analysis. Wiley Series in Probability
  and Statistics. John Wiley \& Sons, Ltd., Chichester; 2000.

\bibitem{sobol1993}
Sobol~IM. Sensitivity estimates for nonlinear mathematical models. Mathematical
  modelling and computational experiments. 1993;\hspace{0pt}1(4):407--414.

\bibitem{gamboa2013_s}
Gamboa~F, Janon~A, Klein~T, Lagnoux~A, Prieur~C. Statistical inference for
  $\textrm{Sobol}$ pick-freeze $\textrm{Monte Carlo}$ method. Statistics.
  2016;\hspace{0pt}50(4):881--902.

\bibitem{rugama2018}
Rugama~LAJ, Gilquin~L. Reliable error estimation for $\textrm{Sobol}´$
  indices. Statistics and Computing. 2018;\hspace{0pt}28(4):725--738.

\bibitem{klein2014}
Gamboa~F, Janon~A, Klein~T, Lagnoux~A. Sensitivity analysis for
  multidimensional and functional outputs. Electron J Statist.
  2014;\hspace{0pt}8(1):575--603;
  \urlprefix\url{https://doi.org/10.1214/14-EJS895}.

\bibitem{marrel2011}
Marrel~A, Iooss~B, Jullien~M, Laurent~B, Volkova~E. Global sensitivity analysis
  for models with spatially dependent outputs. Environmetrics.
  2011;\hspace{0pt}22(3):383--397.

\bibitem{daveiga2015}
Da~Veiga~S. Global sensitivity analysis with dependence measures. Journal of
  Statistical Computation and Simulation. 2015;\hspace{0pt}85(7):1283--1305.

\bibitem{pianosi2015}
Pianosi~F, Wagener~T. A simple and efficient method for global sensitivity
  analysis based on cumulative distribution functions. Environmental Modelling
  \& Software. 2015;\hspace{0pt}67:1--11.

\bibitem{liu2006r}
Liu~H, Chen~W, Sudjianto~A. Relative entropy based method for probabilistic
  sensitivity analysis in engineering design. Journal of Mechanical Design.
  2006;\hspace{0pt}128(2):326--336.

\bibitem{borgonovo2007}
Borgonovo~E. A new uncertainty importance measure. Reliability Engineering \&
  System Safety. 2007;\hspace{0pt}92(6):771--784.

\bibitem{borgonovo2017}
Borgonovo~E. Sensitivity analysis: An introduction for the management
  scientist. Vol. 251. Springer; 2017.

\bibitem{borgonovo2016}
Borgonovo~E, Hazen~GB, Plischke~E. A common rationale for global sensitivity
  measures and their estimation. Risk Analysis.
  2016;\hspace{0pt}36(10):1871--1895.

\bibitem{gamboa2015}
Gamboa~F, Klein~T, Lagnoux~A. Sensitivity analysis based on
  $\textrm{Cram\'er-von Mises}$ distance. SIAM/ASA Journal on Uncertainty
  Quantification. 2018;\hspace{0pt}6(2):522--548.

\bibitem{rao1945}
Rao~CR. Information and the accuracy attainable in the estimation of
  statistical parameters. Bull Calcutta Math. 1945;\hspace{0pt}:81--91.

\bibitem{mardia1972}
Mardia~KV. Statistics of directional data. Academic Press.; 1972.

\bibitem{bhattacharya2012}
Bhattacharya~A, Bhattacharya~R. Nonparametric inference on manifolds: with
  applications to shape spaces. Vol.~2. Cambridge University Press; 2012.

\bibitem{ellingson2015}
Patrangenaru~V, Ellingson~L. Nonparametric statistics on manifolds and their
  applications to object data analysis. CRC Press; 2015.

\bibitem{gamboa2015b}
Chastaing~G, Gamboa~F, Prieur~C. Generalized $\textrm{Sobol}$ sensitivity
  indices for dependent variables: numerical methods. Journal of Statistical
  Computation and Simulation. 2015;\hspace{0pt}85(7):1306--1333.

\bibitem{julia2012}
Bezanson~J, Karpinski~S, Shah~VB, Edelman~A. Julia: A fast dynamic language for
  technical computing. arXiv preprint arXiv:12095145. 2012;\hspace{0pt}.

\bibitem{team2017}
Team~RC. R: A language and environment for statistical computing. Vienna,
  Austria; 2017.

\bibitem{docarmo1992}
do~Carmo~M. Riemannian geometry. {Mathematics. Birkh{\"a}user}, Boston Basel
  Berlin; 1992.

\bibitem{petersen2006}
Petersen~P. Riemannian geometry. Vol. 171. Springer; 2006.

\bibitem{folland2013}
Folland~GB. Real analysis: modern techniques and their applications. John Wiley
  \& Sons; 2013.

\bibitem{pennec2006}
Pennec~X. Intrinsic statistics on $\textrm{Riemannian manifolds}$: Basic tools
  for geometric measurements. Journal of Mathematical Imaging and Vision.
  2006;\hspace{0pt}25(1):127.

\bibitem{christensen1970}
Christensen~JPR. On some measures analogous to haar measure. Mathematica
  Scandinavica. 1970;\hspace{0pt}26(1):103--106.

\bibitem{fraiman2017}
Fraiman~R, Gamboa~F, Moreno~L. Connecting pairwise geodesic spheres by depth:
  $\textrm{DCOPS}$. Journal of Multivariate Analysis. 2019;\hspace{0pt}169:81--
  -- 94;
  \urlprefix\url{http://www.sciencedirect.com/science/article/pii/S0047259X17307340}.

\bibitem{nash1954}
Nash~J. $\textrm{C}^1$--isometric imbedding. Annals of mathematics.
  1954;\hspace{0pt}:383--396.

\bibitem{sobol2001}
Sobol~IM. Global sensitivity indices for nonlinear mathematical models and
  their $\textrm{Monte Carlo}$ estimates. Mathematics and computers in
  simulation. 2001;\hspace{0pt}55(1):271--280.

\bibitem{klein2013}
Janon~A, Klein~T, Lagnoux~A, Nodet~M, Prieur~C. Asymptotic normality and
  efficiency of two $\textrm{Sobol}$ index estimators. ESAIM: Probability and
  Statistics. 2014;\hspace{0pt}18:342--364.

\bibitem{ibragimov1999}
Ibragimov~R, Sharakhmetov~S. {Analogues of $\textrm{Khintchine,
  Marcinkiewicz-Zygmund and Rosenthal}$ Inequalities for Symmetric Statistics}.
  Scandinavian Journal of Statistics. 1999;\hspace{0pt}26(4):621--633.

\bibitem{gine1992b}
Arcones~MA, Gine~E. On the bootstrap of $\textrm{U}$ and $\textrm{V}$
  statistics. The Annals of Statistics. 1992;\hspace{0pt}:655--674.

\bibitem{kendall1974}
Kendall~DG. Pole-seeking brownian motion and bird navigation. Journal of the
  Royal Statistical Society Series B (Methodological).
  1974;\hspace{0pt}:365--417.

\bibitem{watson1983}
Watson~GS. Statistics on spheres. Wiley; 1983.

\bibitem{landau1965}
Landau~LD, Lifshitz~EM. Theory of elasticity. Nauka; 1965.

\bibitem{moakher2005}
Moakher~M. A differential geometric approach to the geometric mean of symmetric
  positive-definite matrices. SIAM Journal on Matrix Analysis and Applications.
  2005;\hspace{0pt}26(3):735--747.

\bibitem{billingsley1999}
Billingsley~P. Convergence of probability measures. John Wiley \& Sons; 2013.

\bibitem{hoeffding1963}
Hoeffding~W. Probability inequalities for sums of bounded random variables.
  Journal of the American Statistical Association.
  1963;\hspace{0pt}58(301):13--30.

\bibitem{serfling1980}
Serfling~RJ. Approximation theorems of mathematical statistics. Vol. 162. John
  Wiley \& Sons; 2009.

\end{thebibliography}

\end{document}